\documentclass[12pt]{amsart}
\usepackage{amscd,amssymb,showkeys}

\usepackage{curves}
\usepackage{float}

\theoremstyle{latex 2e}
\newtheorem{thm}[subsection]{Theorem}
\newtheorem{lem}[subsection]{Lemma}
\newtheorem{prop}[subsection]{Proposition}
\newtheorem{cor}[subsection]{Corollary}

\newtheorem{defn}[subsection]{Definition}

\theoremstyle{remark}
\newtheorem*{rem}{Remark}

\numberwithin{equation}{section}

\begin{document}

\today

\title[Perelman's $\lambda$-functional and the Seiberg-Witten equations] %
{Perelman's $\lambda$-functional and the Seiberg-Witten equations}

\author{Fuquan Fang}
\thanks{The first author was supported by
NSF Grant 19925104 of China, 973 project of Foundation Science of
China, and the Capital Normal University}
\address{Nankai Institute of Mathematics,
Weijin Road 94, Tianjin 300071, P.R.China}
\address{Department of Mathematics, Capital Normal University,
Beijing, P.R.China}
  \email{ffang@nankai.edu.cn}
\author{Yuguang Zhang}
\address{Department of Mathematics, Capital Normal University,
Beijing, P.R.China  }

\begin{abstract}

In this paper we study the supremum of Perelman's
$\lambda$-functional ${\lambda }_M(g)$ on Riemannian $4$-manifold
$M$ by using the Seiberg-Witten equations. We prove among others
that, for a compact K\"{a}hler-Einstein complex surface $(M, J,
g_{0})$ with negative scalar curvature, (i) If $g_{1}$ is a
Riemannian metric on $M$  with $\lambda_{M}(g_{1})=
\lambda_{M}(g_{0})$, then $\text{Vol}_{g_{1}}(M)\geq
\text{Vol}_{g_{0}}(M)$. Moreover, the equality holds if and only
if $g_{1}$ is also a K\"{a}hler-Einstein metric with negative
scalar curvature. (ii) If $\{g_{t}\}$, $t\in [-1,1]$,  is  a
family of Einstein metrics on $M$ with initial metric $g_{0}$,
then $g_{t}$ is a K\"{a}hler-Einstein metric with negative scalar
curvature.
\end{abstract}
\maketitle

\section{Introduction }

In his celebrated paper [H] R.Hamilton introduced the Ricci-flow
evolution equation
\begin{equation}
\frac{\partial}{\partial t}g(t)=-2{\rm Ric}(g(t))
\end{equation}
with initial metric $g(0)=g$. The Ricci flow is now a fundamental
tool to solve the famous Poincar\'e conjecture and Thurston's
Geometrizaion conjecture, by the works of G. Perelman [Pe1][Pe2].
A fundamental new discovery of Perelman is to prove the Ricci-flow
evolution equation is the gradient flow of a so called {\it
Perelman's $\lambda$-functional} of a Riemannian manifold(cf.
[Pe1][KL]), which may be described as follows: for a smooth
function $f\in C^\infty(M)$ on a Riemannian $n$-manifold with a
Riemannian metric $g$, let
\begin{equation}\mathcal{F}(g,f)=\int_{M} (R_g+|\nabla f|^{2})e^{-f}
dvol_{g},\end{equation} where $R_g$ is the scalar curvature of
$g$.
 The Perelman's
$\lambda$-functional  is defined by
\begin{equation}\lambda_{M}(g)=\inf _f \{\mathcal{F}(g,f)|\int_{M}e^{-f}
 dvol_{g}=1\}.\end{equation} Note that $\lambda_{M}(g)$ is the lowest eigenvalue of the
 operator $-4\triangle+R_{g}$.  Let \begin{equation}\overline{\lambda}_{M}(g)
 =\lambda_{M}(g)\text{Vol}_{g}(M)^{\frac{2}{n}}\end{equation}
which is invariant up to rescale the metric.  Perelman [Pe1] has
established the monotonicity property of $\overline{\lambda}
_M(g_t)$ along the Ricci flow $g_t$, namely, the function is
non-decreasing along the Ricci flow $g_t$ whenever
$\overline{\lambda} _M(g_t)\leq 0$. Therefore, it is interesting to
study the upper bound of $\overline{\lambda}_{M}(g)$. This leads to
define a diffeomorphism invariant
  $\overline{\lambda}_{M} $ of $M$ due to Perelman (cf. [Pe2]
  [KL])
by \begin{equation}\overline{\lambda}_{M}=
  \sup\limits_{g\in \mathcal{M}} \overline{\lambda}_{M}(g),
  \end{equation} where $\mathcal{M}$ is the set of Riemannian
metrics on $M$. It is easy to see that $\overline{\lambda} _M=0$
if $M$ admits a volume collapsing with bounded scalar curvature
but does not admit any metric with positive scalar curvature (cf.
[KL]). By a deep result of Perelman (cf. [Pe2][KL]), for a
$3$-manifold $M$ which does not admit a metric of positive scalar
curvature, $(-\overline{\lambda}_{M})^{\frac 32}$ is proportional
to the minimal volume of the manifold.

The invariant $\overline{\lambda}_{M} $ may take value $+\infty$,
e.g., $M=S^2\times S^2$. Thus, it seems only interesting when
$\overline{\lambda}_{M} \leq 0$, i.e, when $M$ does not admit any
metric of positive scalar curvature.

In this paper we will investigate $\overline{\lambda}_{M} $ by
using the Seiberg-Witten monopole equations for a $4$-manifold
$M$. We say a $\rm Spin^{c}$-structure (or equivalently its first
Chern class) is a {\it monopole class} if the Seiberg-Witten
monopople equations has an irreducible solution.

Our first result is as follows:

\begin{thm} Let $(M, \mathfrak{c})$ be a smooth  compact closed
oriented $ 4$-manifold with a  $\rm Spin^{c}$-structure
$\mathfrak{c}$. If the first Chern class $c_1$ of $\mathfrak{c}$
is a monopole class of $M$ satisfying that  $c_{1}^{2}[M]
>0$.  Then, for any Riemannian metric $g$,
\begin{equation} \overline{\lambda}_{M}(g)\leq - \sqrt{32\pi^{2}
c_{1}^{2}[M]}.\end{equation}  Moreover, the equality holds if and
only if $g$ is a K\"{a}hler-Einstein metric with negative scalar
curvature. \end{thm}

\vskip 2mm

By [Ta] the canonical class of a symplectic manifold is a monopole
class. Thus, Theorem 1.1 applies to a K\"ahler minimal surface of
general type, since by [BHPV] $K_X^2>0$ if $X$ is a minimal surface
of general type. We remark that Theorem 1.1 implies that
$\overline{\lambda}_M$ {\it is not a topological invariant of the
underlying manifold.} Indeed, for any pair of positive integers $(m,
n)$, so that $\frac nm \in (\frac 15, 2)$, by [BHPV] VII Theorem 8.3
there is a simply connected minimal surface $X$ of general type so
that $m=c_2(X), n=c_1^2(X)$. Let $M$ be the blow up of $X$ at one
point. Then $c_1^2[M]=n-1$. By Theorem 1.1 we know that
$\overline{\lambda}_{M}(g)\leq - \sqrt{32\pi^{2} (n-1)}$. On the
other hand, since $M$ is a simply connected $4$-manifold of odd
intersection type, by Freedman's classification it is homeomorphic
to the connected sums  $k\Bbb CP^2\# l\overline{\Bbb CP^2}$ for some
positive integers $k, l$. Since the latter admits a metric with
positive scalar curvature, $\overline{\lambda }_{k\Bbb CP^2\# l
\overline{\Bbb CP^2}}>0$. This shows that $\overline{\lambda}_{M}$
is not a topological invariant.

A geometric consequence of Theorem 1.1 is the following
comparison theorem.

\begin{cor}
 Let $(M, J, g_{0})$ be a  compact
K\"{a}hler-Einstein complex surface with negative scalar
curvature. If $g_{1}$ is a Riemannian metric on $M$  with
$\lambda_{M}(g_{1})= \lambda_{M}(g_{0})$, then
$$\text{Vol}_{g_{1}}(M)\geq \text{Vol}_{g_{0}}(M)$$
Moreover, the equality holds if and only if $g_{1}$ is a
K\"{a}hler-Einstein metric with negative scalar curvature.
\end{cor}

One may wonder whether $(M, g_1)$ and $M, g_0)$ are isometric in
the above theorem when the equality holds. This may not be true.
Indeed, there are infinitely many families of K\"ahler Einstein
metrics on a compact complex surface $M$ in different isometry
classes with negative scalar curvature but all the same volume and
same $\lambda _M(\cdot )$.

The following corollary shows a deformation rigidity of Einstein
metrics on compact complex K\"ahler-Einstein surface with negative
scalar curvature.

\begin{cor} Let $(M, J, g_{0})$ be a  compact K\"{a}hler-Einstein complex
surface with negative scalar curvature. If $\{g_{t}\}$, $t\in
[-1,1]$,  is  a family of Einstein metrics on $M$ with initial
metric $g_{0}$, then,  for any  $t$,  $g_{t}$ is a
K\"{a}hler-Einstein metric with negative scalar curvature.
\end{cor}

The above Corollary 1.3 should be compared with Corollary D in
[G], where the same conclusion was obtained when $g_0$ has
positive scalar curvature. On the other hand, under some
additional technical assumptions similar results are obtained in
general dimensions in [DWW] and [Ko] along a completely different
line.

For a compact symplectic $4$-manifold $N$ with first Chern class
$c_1$, the Riemann-Roch formula implies that
$c_1^2[N]=2\chi(N)+3\tau(N)$, where $\chi(N)$ and $\tau(N)$ are the
Euler characteristic and the signature of $N$ respectively. If $N$
admits a  K\"{a}hler-Einstein metric  with negative scalar
curvature, we have already known  $\overline{\lambda}_{N}= -
(32\pi^{2} (2\chi(N)+3\tau(N)))^{\frac{1}{2}}$. In the  next
theorem, we will show the exact quantity of $\overline{\lambda}_{M}$
where $M$
 is obtained by blowing-up $N$ at $k$ points.
\vskip 2mm

\begin{thm} Let $(N, \omega)$ be  a compact symplectic 4-manifold with
$b^{+}_{2}(N)   >1$. Let $M=N\sharp k \overline{\mathbb{C}P^{2}}$,
where $k\ge 0$. If $2\chi(N)+3\tau(N)>0$, then
\begin{equation} \overline{\lambda}_{M}\leq -\sqrt{32\pi^{2}
(2\chi(N)+3\tau(N))}.
  \end{equation} Furthermore, the equality holds if $N$ admits a K\"{a}hler-Einstein metric.
  \end{thm}

\vskip 2mm

 It is known that  the Seiberg-Witten invariant of connected sums vanishes if
 both factors have positive $b_2^+$. In [Ba] [BaF],  a refinement of the Seiberg-Witten
invariant is defined, which may not vanish for connected sums of
few factors. This may be used to improve the above theorem as
follows:

\vskip 2mm

\begin{thm} Let $(N_{i}, \omega_{i})$, $i=1, \cdots ,\ell$, where $\ell \le 4$,
be compact symplectic 4-manifolds satisfying that
$b_{1}(N_{i})=0$, $b^{+}_{2}
   (N_{i})\equiv 3 ( {\rm mod}4)$, and $\sum_{i=1}^{4}b^{+}_{2}
   (N_{i})\equiv 4( {\rm mod}8)$. Assume that $c_{1}^{2}[N_{1}] >0$, and
   $c_{1}^{2}[N_{i}] \geq 0$ for all $i$. Let $X$ be a compact oriented
    4-manifold with $b^{+}_{2}(X)=0$, which admits a metric of positive scalar
     curvature. Let $M=\sharp_{i=1}^{\ell}N_{i}\sharp X$. Then
    \begin{equation} \overline{\lambda}_{M}\leq -\sqrt{32\pi^{2}
\sum_{i=1}^{\ell}c_{1}^{2}[N_{i}]}
  \end{equation}  Furthermore, the equality holds  if  $N_{1}, \cdots, N_{\ell}$ admit  K\"{a}hler-Einstein
  metrics.
   \end{thm}

\vskip 2mm

The technique developed in proving Theorem 1.4 and 1.5 has an easy
corollary:

\vskip 2mm

\begin{cor}
Let $(N_{i}, g_{i})$, $i=1,\cdots, l_{1}$,
   be  compact Riemannian 4m-manifolds $(m\geq 2)$  with holonomy $SU(2m)$ or $Sp(m)$
  or $Spin(7)$, and  $X_{j}$, $j=1,\cdots, l_{2}$, be
simply connected  compact oriented spin 4m-manifolds  with vanishing
$\widehat{A}$-genus,   $\widehat{A}(X_{i})=0$.
   If $M=\sharp_{i=1}^{l_{1}}N_{i}\sharp \sharp_{j=1}^{l_{2}} X_{j}
    $, and $\widehat{A}(M)\neq 0$,   then $$ \overline{\lambda}_{M}=0.$$

\end{cor}

\vskip 2mm

Let $M$ be a smooth  compact oriented $4$-manifold. By Perelman
[Pe1] a critical point of $\overline{\lambda}_M(\cdot )$ is an
Einstein metric. Therefore, it is interesting to ask

\noindent {\bf Question}: {\it Can one deform a metric $g$ to an
Einstein metric through the Ricci flow, provided
$\overline{\lambda}_M(g)$ is sufficiently close to the maximum
$\overline{\lambda _M}$ of the $\lambda$-functional?}

This may not have a positive answer in general, of course,  e.g.,
for a graph $3$-manifold $M$, by [Pe2][KL] $\overline{\lambda
_M}=0$, but $M$ can not have any Einstein metric except $M$ is a
flat manifold.

To formulate our next result, let us consider the moduli space of
metrics
   $$\mathcal{M}_{(\Lambda, D)}=\{g:  |K_{g}|<\Lambda^{2}, diam_{g}< D\}, $$ where  $diam_{g}$ is the diameter, and
      $K_{g}$ is the sectional curvature of $g$.

\vskip 2mm

\begin{prop} Let $(M, J)$ be a  compact
 almost complex  $4$-manifold satisfying that $\chi (M) \in [\frac 32 \tau (M),  3\tau (M)]$, and
 $\tau (M)>0$. If
  the canonical $\rm Spin^{c}$-structure
$\mathfrak{c}$ induced by $J$ is a monopole class, then there exists
a constant $\varepsilon=\varepsilon(\Lambda, D)>0$ depending only on
$\Lambda$ and   D  such that for any Riemannian metric $g\in
\mathcal{M}_{(\Lambda, D)}$ on $M$ satisfying that
\begin{equation}\overline{\lambda}_{M}(g)\geq - \sqrt{32\pi^{2}
(2\chi(M)+3\tau (M) )}-\varepsilon \end{equation} it can be
deformed to a complex  hyperbolic metric  through the Ricci flow.
\end{prop}

\vskip 2mm

The rest of the paper is organized as follows:  In $\S2$ we recall
some facts about Seiberg-Witten equations. In $\S3$ we prove Theorem
1.1, Corollary 1.2 and Corollary 1.3. In $\S4$ we prove Theorem 1.4,
 Theorem  1.5 and Corollary  1.6. In $\S5$ we prove Proposition  1.7.

\vskip 8mm

\noindent {\bf Acknowledgement:}
    The second author is grateful to Zhenlei
  Zhang for explaining Perelman's work to him.

\vskip 8mm

\section{Preliminaries }

 In this section, we recall some facts about Seiberg-Witten
 equations. More details can be found in [N1] and [Le2].

Let $(M, g)$  be a   compact oriented  Riemannian  $4$-manifold
with a  $\rm Spin^{c}$ structure $\mathfrak{c}$. Let $b^{+}_{2}$
denote the dimension of the space of self-dual harmonic $2$-forms
in $M$. Let $S^{\pm}_{\mathfrak{c}}$ denote the $\rm
Spin^{c}$-bundles associated to $\mathfrak{c}$, and let
 $L$ be the determinant line bundle of $\mathfrak{c}$. There is a well-defined
 Dirac operator $$\mathcal{D}_{A}:
  \Gamma(S^{+}_{\mathfrak{c}})\longrightarrow
  \Gamma(S^{-}_{\mathfrak{c}})$$

Let $c:
  \wedge^{*}T^{*}M \longrightarrow {\rm End}(S^{+}_{\mathfrak{c}}\oplus
  S^{-}_{\mathfrak{c}})$ denote the Clifford multiplication on the $\rm{Spin}^c$-bundles,
and, for any $\phi\in \Gamma(S^{\pm})$, let
  $$q(\phi)=\overline{\phi}\otimes\phi-\frac{1}{2}|\phi|^{2}{\rm id}.$$
   The Seiberg-Witten equations read
   \begin{equation}\begin{array}{ccc}\mathcal{D}_{A}\phi=0 \\
   c(F^{+}_{A})=q(\phi)
\end{array}   \end{equation}
 where the unknowns are  a hermitian
 connection $A$ on $L$ and a section $\phi\in
 \Gamma(S^{+}_{\mathfrak{c}})$., and $F^{+}_{A}$ is the self-dual
 part of the curvature of $A$.

 A resolution of (2.1) is called  {\it reducible} if $\phi\equiv 0$;
 otherwise, it is called {\it irreducible}. If $(\phi, A)$ is a resolution of (2.1),
  then  one calculates   \begin{equation}|F^{+}_{A}|=\frac{1}{2\sqrt{2}}
 |\phi|^{2},  \end{equation}

The Bochner formula reads
\begin{equation} 0=2\triangle |\phi|^{2}
 +4|\nabla^{A}\phi|^{2}+R_{g}|\phi|^{2}+|\phi|^{4}, \end{equation} where
 $R_{g}$ is the scalar curvature of $g$.

The Seiberg-Witten invariant can be defined by counting the
irreducible solutions of the Seiberg-Witten equations (cf. [N1]
[Le2]).

\begin{defn}[K1] Let $M$ be a smooth compact
  oriented   $4$-manifold. An element
  $\alpha\in H^{2}(M, \mathbb{Z})/$torsion is called a monopole
  class of $M$ if and only if there exists a $\rm Spin^{c}$-structure
  $\mathfrak{c}$ on $M$ with the first Chern class $c_{1}\equiv
  \alpha$\rm(mod torsion), \it so that the Seiberg-Witten equations have a solution for every
   Riemannian metric $g$ on $M$.
\end{defn}

Deep results have been found in Seiberg-Witten theory to detect
the monopole classes. For example, if  $(M, \omega)$ is a compact
symplectic 4-manifold with $b^{+}_{2}>1$,  the canonical class of
$(M, \omega)$ is a monopole class (cf. [Ta], or Theorem 4.2 in
   [K2]).

    A refinement of Seiberg-Witten
 invariant is defined in [Ba] [BaF], which takes values in a
 cohomotopy group. The remarkable fact is that this invariant is not
 killed off by the sort of connected sum operation.  If  $(N_{i}, \omega_{i})$, $i\in
\{1,2,3,4\}$,  are  the same as in Theorem 1.5, then,  by
Proposition 10 in [IL],  $\sum_{i=1}^{l}\pm c_{1}(N_{i})$ is a
monopole
  class of $\sharp_{i=1}^{l}N_{i}$ where  $l\leq 4$.

\vskip 5mm

\section{Proof of Theorems 1.1}

Let $(M, g)$ be a  Riemannian  $\rm Spin^{c}$-manifold of
dimension $n$. To  prove  Theorem 1.1, we need the following
version of Kato's inequality.

\begin{lem} Let $\phi$ be a harmonic $\rm Spin^{c}$-spinor on $(M, g)$, i.e. $\mathcal{D}_{A}\phi=0$, where
 $\mathcal{D}_{A}$ is the Dirac operator and $A$ is a connection on
 the determinant line bundle. Then \begin{equation}|\nabla|\phi||^{2}\leq |\nabla^{A}\phi|^{2}
 \end{equation} at all points where $\phi$ is non-zero. Moreover,
  the equality can only occur if $\nabla^{A}\phi\equiv 0$.
\end{lem}

\begin{proof} Fix a point $p\in M$ at which $\phi(p)\neq 0$ so
 that $|\phi|$ is differentiable at $p$.  If  $e_{1}, \cdots, e_{n}$
    is  an orthonormal basis of $T_{p}M$, then
    $|\nabla|\phi||^{2}=\sum |\frac{\partial}{\partial
    e_{i}}|\phi||^{2}$ and $|\nabla^{A}\phi|^{2}=\sum
    |\nabla^{A}_{e_{i}}\phi|^{2}$. For any $i$,  $$|\phi||\frac{\partial}{\partial
    e_{i}}|\phi||=\frac{1}{2}|\frac{\partial}{\partial
    e_{i}}|\phi|^{2}|=|{\rm Re} \langle \nabla^{A}_{e_{i}}\phi, \phi
    \rangle|\leq |\nabla^{A}_{e_{i}}\phi||\phi|,$$ $$|\frac{\partial}{\partial
    e_{i}}|\phi||\leq |\nabla^{A}_{e_{i}}\phi|, \ \ \ \ {\rm and} \ \ \  |\nabla|\phi||^{2}\leq |\nabla^{A}\phi|^{2}.$$ The  equality can only occur
    if there are  real  numbers  $\alpha_{i}$ such that
    $\nabla^{A}_{e_{i}}\phi=\alpha_{i} \phi$. Since $\phi$ is a harmonic $\rm Spin^{c}$
    spinor, $$0=\mathcal{D}_{A}\phi=\sum c(e_{i})
    \nabla^{A}_{e_{i}}\phi=c(\sum \alpha_{i}e_{i}) \phi=c(w)
    \phi,$$ where $w=\sum \alpha_{i}e_{i}$ and $c$ is the Clifford multiplication. Then $$0=-|w|^{2}\phi.$$
     Thus $w=0$, $\alpha_{i}=0$ and $\nabla^{A}\phi=0$ at  $p$.
     Thus  we obtain the conclusion.
\end{proof}

      For any $\varepsilon >0$, let
      $|\phi|_{\varepsilon}^{2}=|\phi|^{2}+\varepsilon^{2}$. If
      $\phi$ is harmonic, by above lemma, \begin{equation}
      |\nabla|\phi|_{\varepsilon}|^{2}\leq \frac{|\phi|}{|\phi|_{\varepsilon}}|\nabla|\phi||^{2}
      \leq |\nabla^{A}\phi|^{2} \end{equation} at points where
      $\phi(p)\neq 0$.  Since $\{p\in M| \phi(p)\neq 0\} $ is dense in $M$ for harmonic
      $\phi$, we conclude that (3.2) holds  everywhere in $M$.

\begin{prop} Let $(M, g)$ be a  compact oriented
 Riemannian   $4$-manifold, and $\mathfrak{c}$ be a  $\rm Spin^{c}$-structure on
 $M$. If there is an irreducible  solution $(\phi, A)$ to the Seiberg-Witten
 equations (2.1) for $g$ and $\mathfrak{c}$, then
\begin{equation} \overline{\lambda}_{M}(g)\leq - \sqrt{32\pi^{2}[c_{1}^{+}]^{2}[M]},\end{equation} where
$c_{1}^{+}$ is the self-dual part of  the harmonic form
representing the first Chern class $c_{1}$ of $\mathfrak{c}$. When
$[c_{1}^{+}]\neq 0$, equality can only occur if $g$ is a
K\"{a}hler metric with constant  negative scalar curvature.
\end{prop}

\begin{proof} Let $(\phi, A)$ be an irreducible  solution to the
Seiberg-Witten
 equations.  The Bochner formula implies
$$0=\frac{1}{2}\Delta
|\phi|^{2}+|\nabla^{A}\phi|^{2}+\frac{R_{g}}{4}|\phi|^{2}+\frac{1}{4}|\phi|^{4},$$
Therefore
$$\int_{M}(|\nabla^{A}\phi|^{2}+\frac{R_{g}}{4}|\phi|^{2})dvol_{g}=-\frac{1}{4}
\int_{M}|\phi|^{4}dvol_{g}.$$ By (3.2),
$$\int_{M}(|\nabla|\phi|_{\varepsilon}|^{2}+\frac{R_{g}}{4}|\phi|_{\varepsilon}^{2})dvol_{g}\leq -\frac{1}{4}
\int_{M}|\phi|^{4}dvol_{g}+\varepsilon^{2}\int_{M}\frac{R_{g}}{4}dvol_{g}.$$
Since $\lambda_{M}(g)$ is the lowest eigenvalue of the
 operator $-4\triangle+R_{g}$, we obtain $$
 \lambda_{M}(g)\int_{M}|\phi|_{\varepsilon}^{2}dvol_{g}\leq
  \int_{M}(4|\nabla|\phi|_{\varepsilon}|^{2}+R_{g}|\phi|_{\varepsilon}^{2})dvol_{g}.$$
 Note that, for $\varepsilon\ll 1$, $ \lambda_{M}(g)\leq 0$.   By Schwarz inequality, \begin{eqnarray*}
 \lambda_{M}(g)\text{Vol}_{g}(M)^{\frac{1}{2}}(\int_{M}|\phi|_{\varepsilon}^{4}dvol_{g})^{\frac{1}{2}}
  & \leq & \lambda_{M}(g)\int_{M}|\phi|_{\varepsilon}^{2}dvol_{g}\\ & \leq & -
\int_{M}|\phi|^{4}dvol_{g}+\varepsilon^{2}\int_{M}R_{g}
dvol_{g}.\end{eqnarray*} Letting  $\varepsilon \longrightarrow 0$,
we obtain
$$\overline{\lambda}_{M}(g)=\lambda_{M}(g)\text{Vol}_{g}(M)^{\frac{1}{2}}\leq -(
\int_{M}|\phi|^{4}dvol_{g})^{\frac{1}{2}}.$$ By the second
equation in the Seiberg-Witten equations, we get that
$$\overline{\lambda}_{M}(g)\leq
-(
\int_{M}|\phi|^{4}dvol_{g})^{\frac{1}{2}}=-(8\int_{M}|F^{+}_{A}|^{2}dvol_{g})^{\frac{1}{2}}.$$
Note that  $c_{1}^{+}$ is the self-dual part of the harmonic form
representing the first Chern class $c_{1}$. Clearly
$F^{+}_{A}-2\pi c_{1}^{+}$ is $L^{2}$-orthogonal to the harmonic
forms space. Thus
 $$\int_{M}|F^{+}_{A}|^{2}dvol_{g}\geq
 4\pi^{2}\int_{M}|c_{1}^{+}|^{2}dvol_{g}=
 4\pi^{2}\int_{M}c_{1}^{+}\wedge
 c_{1}^{+}=4\pi^{2}[c_{1}^{+}]^{2}$$
 $$\overline{\lambda}_{M}(g)\leq
-\sqrt{32\pi^{2}[c_{1}^{+}]^{2}[M]}$$ If the equality holds, all
of '$\leq$' above are '='.  By Lemma 3.1 and the Bochner formula,
$\nabla^{A}\phi\equiv 0$, $R_{g}=-|\phi|^{2}={\rm const.}$, and
$$\nabla F^{+}_{A}\equiv 0.$$ Note that $F^{+}_{A}$ is a
non-degenerate 2-form since $\phi \ne 0$. Thus, $g$ is a
K\"{a}hler metric with parallel K\"{a}hler form
$\omega=\sqrt{2}\frac{F^{+}_{A}}{|
 F^{+}_{A}|}$.     The desired result follows.
 \end{proof}

\vskip 4mm

\begin{proof}[Proof of Theorem 1.1]  From the hypothesis, for any
Riemannian metric $g$, there is a  solution $(\phi, A)$ of the
Seiberg-Witten
 equations. Let $c_{1}^{+}$ is the self-dual part of the harmonic form
representing the first Chern class $c_{1}$ of $\mathfrak{c}$.
Since
$$\frac{1}{8}\int_{M}|\phi|^{4}dvol_{g}=\int_{M}|F^{+}_{A}|^{2}dvol_{g}\geq
 4\pi^{2}[c_{1}^{+}]^{2}[M]\geq 4\pi^{2}c_{1}^{2}[M]>0,$$ the solution   $(\phi, A)$ is  irreducible.
   By Proposition 3.2, we have $$ \overline{\lambda}_{M}(g)\leq - \sqrt{
32\pi^{2}[c_{1}^{+}]^{2}[M]}\leq -\sqrt{ 32\pi^{2}c_{1}^{2}[M]}.$$
Moreover, if $[c_{1}^{+}]^2[M]\neq 0$, equality can only occur if
$g$ is a K\"{a}hler metric with constant negative scalar
curvature. If $g$ is a metric such that
 the equality  holds in the above formula, then $g$ is a critical point of
the functional $\overline{\lambda}_{M}(\cdot )$.  By the claim in
$\S$2.3 of [Pe1], $g$ is a gradient soliton, i.e, we have the
following equation
$${\rm Ric}(g)-cg+\nabla \nabla f=0,$$ where $c$ is a constant, $f$ satisfies the equation
$$-4\triangle e^{-\frac{1}{2}f}+R_{g} e^{-\frac{1}{2}f}=
\lambda_{M}(g)e^{-\frac{1}{2}f}.$$ Since $\lambda_{M}(g)$ is the
lowest eigenvalue of the
 operator $-4\triangle+R_{g}$ where  $R_{g}$ is a constant,
  we obtain that    $f$ is a constant, and   $g$  is an Einstein metric.

 Now assume that  $g$ is a K\"{a}hler-Einstein metric with negative scalar curvature.   We can assume that the Ricci form
$\rho=-\omega$ where $\omega$ is the K\"{a}hler form associated to
$g$. It is well known that $\rho$ is self-dual and is the harmonic
representative of $2\pi c_{1}$.  We have
$$(32\pi^{2}[c_{1}^{+}]^{2}[M])^{\frac{1}{2}}=(32\pi^{2}c_{1}^{2}[M])^{\frac{1}{2}}
=4(\frac{1}{2}\int_{M}\omega^{2})^{\frac{1}{2}} =4
\text{Vol}_{g}(M)^{\frac{1}{2}}.$$ Since $\lambda_{M}(g)$ is the
lowest eigenvalue of the
 operator $-4\triangle-4$, $\lambda_{M}(g)=-4$. Thus
 $$\overline{\lambda}_{M}(g)=-4\text{Vol}_{g}(M)^{\frac{1}{2}}
 =-\sqrt{ 32\pi^{2}c_{1}^{2}[M]}.$$
 The desired result follows.   \end{proof}

\vskip 3mm

 \begin{proof}[Proofs  of Corollary 1.2 and Corollary 1.3]  From the
hypothesis, $(M, J)$ is a complex surface of general type with
$$c_{1}^{2}[M]=2\chi(M)+3\tau(M)
>0,$$ (cf. Corollary 3.5 in [Le2]).  By Theorem 4.1 in [Le2], the (mod
2) Seiberg-Witten invariant $n_{\mathfrak{c}}(M)\neq 0$ where
$\mathfrak{c}$ is  the  canonical  $\rm Spin^{c}$ structure
induced by $J$.

If $g_{1}$ is a Riemannian metric on $M$, then, by Theorem 1.1,
$$\lambda_{M}(g_{1})\text{Vol}_{g_{1}}(M)^{\frac{1}{2}}
=\overline{\lambda}_{M}(g_{1})\leq - \sqrt{
32\pi^{2}c_{1}^{2}[M]}=\lambda_{M}(g_{0})
\text{Vol}_{g_{0}}(M)^{\frac{1}{2}}.$$ Thus, if
$\lambda_{M}(g_{1})= \lambda_{M}(g_{0})$, we obtain
$$\text{Vol}_{g_{1}}(M)\geq \text{Vol}_{g_{0}}(M)$$ with equality if and only if
$g_{1}$ is a K\"{a}hler-Einstein metric with negative scalar
curvature. This proves Corollary 1.2.

To prove Corollary 1.3, let  $\{g_{t}\}$, $t\in [0,1]$, be  a
family of Einstein metrics starting at $g_{0}$ on $M$, i.e.  ${\rm
Ric}(g_{t})=\frac{R_{g_{t}}}{4}g_{t}$. Let $f_{t}\in
C^{\infty}(M)$ such that, for any $t$,  $e^{-\frac{f_{t}}{2}}$ is
the eigenfunction of the lowest eigenvalue of the operator
$-4\triangle +R_{g_{t}}$ normalized by
$\int_{M}e^{-f_{t}}dvol_{g_{t}}= 1$. Note that
$\lambda_{M}(g_{t})=R_{g_{t}}$, and $f_{t}$ is a constant function
for any $t\in [0,1]$.
 For  a $t_{0}\in [0,1]$, if $v_{ij}=\frac{d}{dt}g_{t,ij}|_{t=t_{0}}$,
 $h=\frac{d}{dt}f_{t}|_{t=t_{0}}$, then we have
 $\frac{d}{dt}dvol_{g_{t}}|_{t=t_{0}}=\frac{1}{2}v dvol_{g_{t_{0}}}$
  where $v=g_{t_{0}}^{ij}v_{ij}$, and $$\int_{M}e^{-f_{t_{0}}}(\frac{1}{2}v-h)dvol_{g_{t_{0}}}=
  0.$$ By the first formula in Section 1 of [Pe1],
  \begin{eqnarray*}
\frac{d}{dt}\lambda_{M}(g_{t})|_{t=t_{0}}& = &
\int_{M}e^{-f_{t_{0}}}[-v_{ij}{\rm
Ric}_{t_{0},ij}+(\frac{1}{2}v-h)R_{g_{t}}]dvol_{g_{t_{0}}}\\
 & = & -\int_{M}e^{-f_{t_{0}}}v_{ij}{\rm
Ric}_{t_{0},ij}dvol_{g_{t_{0}}}\\ & = & -\int_{M}e^{-f_{t_{0}}}v
\frac{R_{g_{t_{0}}}}{4} dvol_{g_{t_{0}}}.
\end{eqnarray*}

\begin{eqnarray*}
\frac{d}{dt}\overline{\lambda}_{M}(g_{t})|_{t=t_{0}}& =
&\frac{d}{dt}\lambda_{M}(g_{t})\text{Vol}_{g_{t_{0}}}(M)^{\frac{1}{2}}+
\frac{1}{2}\lambda_{M}(g_{t_{0}})\text{Vol}_{g_{t_{0}}}(M)^{-\frac{1}{2}}\frac{d}{dt}\text{Vol}_{g_{t}}(M)|_{t=t_{0}}\\
&=&-\text{Vol}_{g_{t_{0}}}(M)^{-\frac{1}{2}}\frac{R_{g_{t_{0}}}}{4}\int_{M}v
dvol_{g_{t_{0}}} \\ & & +
\frac{1}{4}\lambda_{M}(g_{t_{0}})\text{Vol}_{g_{t_{0}}}(M)^{-\frac{1}{2}}\int_{M}v
dvol_{g_{t_{0}}} \\ &= & 0.
\end{eqnarray*}
Hence $$\frac{d}{dt}\overline{\lambda}_{M}(g_{t})\equiv 0.$$ and
so
$$\overline{\lambda}_{M}(g_{t}) \equiv
\overline{\lambda}_{M}(g_{0})=- \sqrt{ 32\pi^{2}c_{1}^{2}[M]},$$
for any  $t$. Therefore, by Theorem 1.1  $g_{t}$ is a
K\"{a}hler-Einstein metric with negative scalar curvature.
Corollary 1.3 follows.
\end{proof}

\vskip 4mm
 \vskip 4mm

\section{Proofs  of Theorem 1.4 and 1.5}

\vskip 3mm

To prove Theorem 1.4 and Theorem 1.5 we need the following
proposition.

\vskip 4mm

\begin{prop} Let $N$ and $X$ be two   smooth  compact
  oriented   n-manifolds, $n\geq 3$, and $M$ be the connected sum of $N$ and $X$, i.e.
  $M=N\sharp X$.\begin{enumerate}\item If $X$  admits a metric with positive scalar
  curvature,  then \begin{equation}\overline{\lambda}_{N}\leq \overline{\lambda}_{M}.
  \end{equation}\item  If $\overline{\lambda}_{N}\leq 0$, $\overline{\lambda}_{X}\leq
  0$,
     $\overline{\lambda}_{M}\leq 0$, and  $n=4m$,   then \begin{equation}
  -(\overline{\lambda}_{X}^{2m}+\overline{\lambda}_{N}^{2m})^{\frac{1}{2m}}\leq \overline{\lambda}_{M}.
  \end{equation}\end{enumerate} \end{prop}

  \vskip 4mm

We remark that the inequality above is often a strict inequality,
e.g. if $N$ is a simply connected Spin-manifold of dimension
$4m\ge 5$ with $\hat A$-genus nonzero and $X=\Bbb CP^{2m}$,
clearly $\overline{\lambda}_{N}\leq 0$, however, by [GL][SY] it is
well-known that $N\#\Bbb CP^{2m}$ admits a metric with positive
scalar curvature, therefore $\overline{\lambda}_{N\# \Bbb
CP^{2m}}>0$.

\vskip 4mm

\begin{lem} Let $(X, h)$ be an oriented compact Riemannian
n-manifold with positive scalar curvature, $N$ be an oriented
smooth compact  n-manifold, $n\geq 3$,  and $M=N\sharp X$. Then,
for any metric $g$ on $N$ and $0< \varepsilon \ll 1$, there exists
a metric $g_{\varepsilon}$ on $M$ such that
$$\lambda_{N}(g)-\varepsilon \leq \lambda_{M}(g_{\varepsilon}),  \ \
\ \ {\rm and} \ \ \ \ \ | \text{Vol}_{g}(N)-
\text{Vol}_{g_{\varepsilon}}(M)|\leq \varepsilon.$$
\end{lem}

\vskip 4mm
\begin{rem} The fact  that $M$ admits a  metric $g_{\varepsilon}$ such that
$\lambda_{M}(g_{\varepsilon})$ is close to $\lambda_{N}(g)$ is an
easy consequence of Theorem 3.1 in [BD2].
 Here  we must construct $g_{\varepsilon}$ carefully such that
 $\text{Vol}_{g_{\varepsilon}}(M) $ is  close to $\text{Vol}_{g}(N)$. \end{rem}

\vskip 4mm

\begin{proof} For a  $p\in N$, denote  $U(r)=\{x| dist_{g}(x, p)<
r\}$. By Lemma 3.7 in [BD2], there exists a $0<\overline{r}< 1$
such that, for any $0<r< \frac{\overline{r}^{11}}{2}$ and any
smooth function $u$ on $ A(r, (2r)^{\frac{1}{11}})$, the following
holds
\begin{equation} \|u\|^{2}_{L^{2}(A(r, 2r))}\leq 10 r^{\frac{5}{2}} \|u\|^{2}_{L^{2}(A(r, (2r)^{\frac{1}{11}}))}
\end{equation} if $\int_{\partial U(\rho)}u \partial_{\nu} u
dA\geq 0$ holds for all $\rho \in [r, (2r)^{\frac{1}{11}}]$. Here
$ A(r, (2r)^{\frac{1}{11}})=\{x| r \leq dist_{g}(x, p)\leq
(2r)^{\frac{1}{11}}\}$, and   $\nu$ is  the unite normal vector
field of $\partial U(\rho)$ pointing away from $p$.  Let
$\Lambda$ be a positive constant bigger than  the lowest
eigenvalue of the operator $-4\triangle + R_{g}$ on $(N\backslash
U(\overline{r}), g)$ with Dirichlet boundary conditions. Let
$R_{0}$ be a lower bound of the scalar curvature $ R_{g}$ of $(N,
g)$, and $R_{1}$ be a number such that
\begin{equation}R_{1}>\min \{0, \Lambda\}, \ \ \ \ \Lambda
\frac{\Lambda-R_{0}}{R_{1}-\Lambda}\leq \frac{\varepsilon}{2}.
\end{equation} By the
arguments in the proof of Theorem 3.1 in [BD2] or Proposition 2.1
of [BD1], there is a metric $g'$ on $N$ arbitrarily close to $g$
in the $C^{1}$-topology such that $R_{g'}\geq R_{0}$ and
$R_{g'}\geq 2R_{1}$ on a neighborhood $U_{0}$ of $p$. Since both
$\lambda_{N}(g)$ and $Vol_{g}(N)$  depend  continuously on $g$ in
the $C^{1}$-topology (See  Lemma 3.4 in [BD2]), we may without
loss of generality  assume that $R_{g}\geq R_{0}$ and $R_{g}\geq
2R_{1}$ on a neighborhood $U_{0}$ of $p$.

Now we choose $r > 0$ and $\zeta >0$ so small that
\begin{enumerate} \item \begin{equation} \frac{R_{1}-R_{0}}{R_{1}-\Lambda} ((\Lambda
+ 1-R_{0})\zeta +\zeta^{2})\leq
\frac{\varepsilon}{2},\end{equation} \item
$2\sqrt{10}r^{\frac{1}{4}}<\zeta$, \item
$U((2r)^{\frac{1}{11}})\subset U_{0}$, \item
$(2r)^{\frac{1}{11}}\leq \overline{r}$, \item
$\text{Vol}_{g}(U(r))<\frac{\varepsilon}{2}$.
\end{enumerate} Let $\eta$ be a smooth cut-off function such that \begin{enumerate}
\item $0\leq \eta \leq 1$ on $N$, \item $\eta\equiv 0$ on $U(r)$,
\item  $\eta\equiv 1$ on $N \backslash U(2r)$, \item $|d \eta|\leq
\frac{2}{r}$ on $N$.
\end{enumerate}
\vskip 4mm

\begin{lem} For any $0<\theta_{0}\ll 1$,  there is a metric
$\widetilde{g}_{\theta_{0}}$ on $A(r, \frac{r}{2})=U(r)\backslash
U(\frac{r}{2})$  satisfying  that
$R_{\widetilde{g}_{\theta_{0}}}\geq R_{1}$,
 $\widetilde{g}_{\theta_{0}}$ agrees with $g$ near the boundary
$\partial U(r)$, and $\widetilde{g}_{\theta_{0}}$ agrees with
$dt^{2}+\delta^{2}g_{0,
 1}$
 near the boundary  $\partial U(\frac{r}{2})\simeq S^{n-1}(1)$,
  where $\delta=\delta(\theta_{0}
)$ is   a function of $\theta_{0}$
 such that $\delta \ll \theta_{0}$, and $g_{0,
 1}$ is  the standard metric of sectional curvature 1
  on $S^{n-1}(1)$. Furthermore,  $$
|\text{Vol}_{\widetilde{g}_{\theta_{0}}}(A(r,
\frac{r}{2}))-\text{Vol}_{g}(U(r))|\longrightarrow 0,$$ if
$\theta_{0} \longrightarrow 0$.   \end{lem} \vskip 4mm

\begin{proof}   We will use Gromov-Lawson's construction here (See
Theorem A of [GL], and  Theorem 3.1 of  [RS]). The key idea of the
proof of Theorem A in  [GL] is to choose a suitable curve $\gamma$
in the t-$\varrho$ plane, and to consider $$T_{\gamma}=\{(t, x)\in
\mathbb{R}\times U(r)|(t, dist_{g}(x, p)=\varrho)\in\gamma  \},$$
with the induced metric, where $\mathbb{R}$ is given the Euclidean
metric and $\mathbb{R}\times U(r)$ is given the natural product
metric $dt^{2}+g$.
 The scalar curvature is given by
\begin{eqnarray*}
R_{\gamma}& =&
R_{g}+((n-1)(n-2)\frac{1}{\varrho^{2}}+O(1))\sin^{2}\theta-(n-1)(\frac{1}{\varrho}+O(\varrho))k\sin
\theta \\ &\geq &
2R_{1}+((n-1)(n-2)\frac{1}{\varrho^{2}}-C)\sin^{2}\theta-(n-1)(\frac{1}{\varrho}+C'\varrho)k\sin
\theta , \end{eqnarray*} where $C$ $C'$ are constants depending
only on the curvature of $g$,  $k$ is the curvature of $\gamma$,
and $\theta$ is the angle between $\gamma$ and the $\varrho$-axis
(See the formula (1) in [GL]). There are several steps to
construct $\gamma$.

\begin{figure}[H]
\setlength{\unitlength}{1cm}
\begin{picture}(14,8)\put(1,2){\makebox(0,0)[t]{$\varrho_{0}$}}
\put(1,1.5){\makebox(0,0)[t]{$\delta$}}\put(3,3){\makebox(0,0)[t]{$\gamma$}}\put(1.75,4.5){\makebox(0,0)[t]{$\theta_{0}$}}\put(1,8){\makebox(0,0)[t]{$\varrho$}}
\put(11,0.75){\makebox(0,0)[t]{$t$}}\put(8,0.75){\makebox(0,0)[t]{$1$}}
\put(4,0.75){\makebox(0,0)[t]{$B$}}\put(1,5){\makebox(0,0)[t]{$\frac{r}{2}$}}
\put(1,1.5){\line
(6,0){10}}\put(1,6){\line(1,-2){3}}\put(1,1){\vector(1,0){10.5}}
\put(1.5,1){\vector(0,1){6.5}}\linethickness{1mm}\qbezier(1.5,8)(1.5,5)(2.5,3)
\qbezier(3,2)(3.25,1.5)(4,1.5)\qbezier(2.5,3)(2.75,2.5)(3,2)\qbezier(4,1.5)(5,1.5)(8,1.5)\end{picture}
\caption{ }
\end{figure}

First, let $\gamma_{0}$ be the bent line segment given by $\{(t,
\varrho)| \varrho=-\coth\theta_{0}t+\frac{r}{2} \}$ on
$\mathbb{R}\times [\frac{r}{4}, 0]$, and a smooth cure with  angle
between $\gamma_{0}$ and the $\varrho$-axis less than $\theta_{0}$
on $\mathbb{R}\times [r, \frac{r}{4}]$. From the proof of Theorem
A in [GL] or the proof of Theorem 3.1 in  [RS], we can choose
$0<\theta_{0}\ll 1$ such that $R_{\gamma_{0}}\geq R_{1}$.

Following  the arguments in P359 of [RS], we choose a
$\varrho_{0}$ with  $0< \varrho_{0}< \min ( \sqrt{\frac{1}{4C}},
\sqrt{\frac{1}{2C'}})$. Then, for $0< \varrho \leq \varrho_{0}$,
we have $$R_{\gamma}\geq  2R_{1}+
(n-1)\frac{3}{4\varrho^{2}}\sin^{2}\theta
-(n-1)\frac{3}{2\varrho}k \sin \theta.$$ Let $\gamma$ be
$\gamma_{0}$ on $\mathbb{R}\times [r,
 \varrho_{0}]$, and be a curve satisfying $k=\frac{\sin
 \theta}{2\varrho}$ on $\mathbb{R}\times [
 \varrho_{0}, 0]$. By the arguments in P359 of  [RS], $\gamma$ is given by the graph of
 function $\varrho=f(t)$ with $f(t)=\delta +
 \frac{1}{4\delta}(t-B)^{2}$.  Note that $(\varrho_{0}, t_{0})\in
 \gamma$, where $t_{0}=(\frac{r}{2}-\varrho_{0})\tan \theta_{0}$,
 and $$\varrho_{0}=\delta +
 \frac{1}{4\delta}(t_{0}-B)^{2}, \ \ \
 f'(t_{0})=\frac{1}{2\delta}(t_{0}-B)=\frac{\varrho_{0}-\frac{r}{2}}{t_{0}}.$$
 Thus we have \begin{equation} \delta=\frac{t^{2}_{0}\varrho_{0}}{(\varrho_{0}-\frac{r}{2})^{2}+t^{2}_{0}},  \ \ \ \
 {\rm and} \ \ \ \ B=t_{0}+\sqrt{4\delta (\varrho_{0}-\delta)}.
 \end{equation} By taking $\theta_{0}\ll r$ and $\varrho_{0}\ll r$,
 we obtain \begin{equation}\delta < 4 \theta_{0}^{2} \varrho_{0}, \ \ \ \
 {\rm and} \ \ \ \ B <r\theta_{0}+4 \theta_{0}
 \varrho_{0}.\end{equation} After $\gamma$ reach $(B, \delta)$,
 let $\gamma$ be $[B, 2B]\times \{\delta\}$. Now we have constructed
 a metric on $T_{\gamma}$, denoted by $g_{\gamma}$, satisfying that $R_{\gamma}\geq  R_{1}$,
   $g_{\gamma}$ agrees with $g$ near $\partial U(r)$, $g_{\gamma}$ agrees
   with the product metric induced by  $\mathbb{R}\times U(r)$ near
   the other
 boundary  of $T_{\gamma}$,
 $\{2B\}\times \partial U(\delta)$. Furthermore, if we let
 $\theta_{0}\longrightarrow 0$, then, by (4.7),
 \begin{equation}|\text{Vol}_{g_{\gamma}}(T_{\gamma})-\text{Vol}_{g}(U(r))|\longrightarrow 0.\end{equation}

 Note that $\partial U(\delta)\cong S^{n-1}(\delta)=\{y\in \mathbb{R}^{n}| \|y\|=\delta\}$.
  If $g_{0,
 1}$ is  the standard metric of sectional curvature 1
  on $S^{n-1}(1)$,  then $\frac{1}{\delta^{2}}g|_{\partial U(\delta)}$ converges to  $g_{0,
 1}$  in the $C^{2}$-topology
 by Lemma 1 in [GL], i.e. there is a 2-tensor $\alpha(\delta)$ on
 $S^{n-1}(1)$
 with $\frac{1}{\delta^{2}}g|_{\partial U(\delta)}-g_{0,
 1}=\alpha(\delta)$ and
 $\|\alpha(\delta)\|_{C^{2}}\longrightarrow 0$ when
 $\delta \longrightarrow 0$. Let $\sigma(t)$ be a smooth  function
 such that $\sigma(t)\equiv 1$ on $[0, \frac{1}{3}]$, $\sigma(t)\equiv 0$ on $[ \frac{2}{3},
 1]$, and $|\frac{d}{dt}\sigma(t)|\leq 4$ on $[0, 1]$.  Define a
 metric on $[0, \frac{1}{\delta}]\times S^{n-1}(1)$ by $g_{\delta}'=dt^{2}+ g_{0,
 1} + \sigma(\delta t)\alpha(\delta)$. When $\delta \ll  1$, $R_{g_{\delta}'}>
 \frac{1}{4}$. Define  $g_{\delta}=\delta^{2}g_{\delta}'$ on $[2B, 1]\times \partial
 U(\delta)$ which satisfies   that $g_{\delta}=dt^{2}+g|_{\partial U(\delta)}$
 near $\{2B\} \times \partial
 U(\delta)$,  $g_{\delta}=dt^{2}+\delta^{2}g_{0,
 1}$
 near $\{1\} \times \partial
 U(\delta)$, $R_{g_{\delta}}> \frac{1}{4\delta^{2}}$, and
 $\text{Vol}_{g_{\delta}}([2B, 1]\times \partial
 U(\delta))=O(\delta^{n-1})=O(\theta_{0}^{2n-2})$. Let $\widetilde{T}_{\gamma}$ be
 the manifold obtained by gluing $T_{\gamma}$ and  $[2B, 1]\times \partial
 U(\delta) $ at $\{2B\}\times \partial
 U(\delta)$,  i.e.
   \begin{equation}\widetilde{T}_{\gamma}=T_{\gamma}\bigcup
  [2B, 1]\times \partial
 U(\delta) , \end{equation} and $\widetilde{g}_{\gamma}$ be a metric
 on $\widetilde{T}_{\gamma}$ such that
 $\widetilde{g}_{\gamma}=g_{\gamma}$ on $T_{\gamma}$, and
 $\widetilde{g}_{\gamma}=g_{\delta}$ on $[2B, 1]\times \partial
 U(\delta)$. Thus the metric $\widetilde{g}_{\gamma}$ satisfies that
 $R_{\widetilde{g}_{\gamma}}\geq R_{1}$ and
  \begin{equation}|\text{Vol}_{\widetilde{g}_{\gamma}}(\widetilde{T}_{\gamma})-\text{Vol}_{g}(U(r))|\longrightarrow 0,\end{equation}
  when $\theta_{0}\longrightarrow 0$. Since $A(r,
\frac{r}{2})\simeq \widetilde{T}_{\gamma}$, we obtain the
conclusion by letting
$\widetilde{g}_{\theta_{0}}=\widetilde{g}_{\gamma}$.
\end{proof} \vskip 3mm

  Let's continue  to prove Lemma 4.2.  Let $\widetilde{U}$ be the connected sum of $U(r)$
and $X$.  Now let's  consider $(X, h)$. By the proof of Theorem A
in [GL],
  we have a compact manifold $\widetilde{X}$ with boundary $\partial
  \widetilde{X}=S^{n-1}(\varsigma)$, which is obtained by deleting a
  small disc from $X$, and a metric $\widetilde{h}$ on
  $\widetilde{X}$ such that the scalar curvature  $R_{\widetilde{h}}$
  is positive, and $\widetilde{h}=dt^{2}+ g_{0, \varsigma}$ near the
  boundary $\partial \widetilde{X}$, where $g_{0,
 \varsigma}$ is  the standard metric of sectional curvature
 $\frac{1}{\varsigma^{2}}$ on $S^{n-1}(\varsigma)$.  By letting  $\theta_{0}\ll \min \{\varsigma, \min
 R_{\widetilde{h}}\}$, we obtain that  $\delta\ll \min \{\varsigma, \min
 R_{\widetilde{h}}\}$, and  the metric
 $(\frac{\delta}{\varsigma})^{2}\widetilde{h}$ satisfies that the
 scalar curvature of $(\frac{\delta}{\varsigma})^{2}\widetilde{h}$
 is bigger than $R_{1}$,
 $(\frac{\delta}{\varsigma})^{2}\widetilde{h}=dt^{2}+ g_{0, \delta}$
  near the
  boundary $\partial \widetilde{X}$, and
  $$\text{Vol}_{(\frac{\delta}{\varsigma})^{2}\widetilde{h}}(\widetilde{X})\longrightarrow
  0,$$ if $\delta\longrightarrow 0$.

  Note that $\widetilde{U}$ is obtained by gluing $A(r, \frac{r}{2})$ and $\widetilde{X}$ at $\partial
 U(\frac{r}{2})\cong \partial\widetilde{X}$,  i.e.  $\widetilde{U}=A(r, \frac{r}{2})\bigcup
  \widetilde{X}$. For any $\theta_{0}\ll 1$, define a metric $\widetilde{g}_{\theta_{0}}'$ on
  $\widetilde{U}$ such that
  $\widetilde{g}_{\theta_{0}}'=\widetilde{g}_{\theta_{0}}$ on
  $A(r, \frac{r}{2})$, where $\widetilde{g}_{\theta_{0}}$ is the metric obtained in
   Lemma 4.3, and $\widetilde{g}_{\theta_{0}}'=(\frac{\delta}{\varsigma})^{2}\widetilde{h}$
  on $\widetilde{X}$, which satisfies that
  $R_{\widetilde{g}_{\theta_{0}}'}\geq R_{1}$ and
  $$|\text{Vol}_{\widetilde{g}_{\theta_{0}}'}(\widetilde{U})-\text{Vol}_{g}(U(r))|
  <|\text{Vol}_{\widetilde{g}_{\theta_{0}}}
  (A(r, \frac{r}{2}))-\text{Vol}_{g}(U(r))|
  + \text{Vol}_{(\frac{\delta}{\varsigma})^{2}\widetilde{h}}(\widetilde{X})\longrightarrow
  0,$$ when $\theta_{0}\longrightarrow 0$.

   Note  that $M$ is obtained by
   gluing $N\backslash U(r)$ and $\widetilde{U}$ at $\partial U(r)$,
   i.e.
   $$M=(N\backslash U(r))\bigcup \widetilde{U}.$$ Define  metrics
   $g_{\theta_{0}}$ on $M$ by $g_{\theta_{0}}=g $ on $N\backslash
   U(r)$ and $g_{\theta_{0}}=\widetilde{g}_{\theta_{0}}'$
   on $\widetilde{U}$, which satisfy  $$|\text{Vol}_{g}(N)-\text{Vol}_{g_{\theta_{0}}}(M)|
   \leq \text{Vol}_{g}(U(r))
   + |\text{Vol}_{g}(U(r))-\text{Vol}_{\widetilde{g}_{\theta_{0}}}(\widetilde{U})| \longrightarrow 0,
   $$ when $\theta_{0}\longrightarrow 0$. Thus, for any $0<
   \varepsilon \ll 1$, there is a $\theta_{0}$ such that \begin{equation}
   |\text{Vol}_{g}(N)-\text{Vol}_{g_{\theta_{0}}}(M)| <\varepsilon.
   \end{equation}  By defining  $g_{\varepsilon}=g_{\theta_{0}}$ on
   $M$, we obtain the volumes inequality.  \vskip 3mm

  \begin{lem} $$\lambda_{N}(g)-\varepsilon \leq
   \lambda_{M}(g_{\varepsilon}).$$ \end{lem}  \vskip 3mm

  \begin{proof} The following arguments is the same as the proof of Theorem 3.1 in [BD2].
   But for reader's convenience, we present the proof here.   Let $u$ be the eigenfunction of
   $\lambda_{M}(g_{\varepsilon})$ on $(M, g_{\varepsilon})$. The
   function $v=\eta u$ can be regarded as a function on $(N, g)$. Thus
\begin{equation}\lambda_{N}(g)\leq \frac{\int_{N}(4|d v|^{2}+R_{g}v^{2})dvol_{g}}{\int_{N}v^{2}dvol_{g}}. \end{equation}
Since  $\Lambda$ is larger than  the lowest eigenvalue of the
operator $-4\triangle + R_{g}$ on $(N\backslash U(\overline{r}),
g)$ with Dirichlet boundary conditions, we have
$\lambda_{M}(g_{\varepsilon})\leq \Lambda$ by Lemma 92.5 in [KL].
Thus
\begin{eqnarray*}R_{1}\int_{\widetilde{U}\bigcup A(r,2r)
}u^{2}dvol_{g_{\varepsilon}}+ R_{0}\int_{M\backslash
\widetilde{U}\bigcup A(r,2r) }u^{2}dvol_{g_{\varepsilon}}& \leq &
\int_{M}(4|d
u|^{2}+R_{g_{\varepsilon}}u^{2})dvol_{g_{\varepsilon}}\\ &\leq &
\Lambda \int_{M}u^{2}dvol_{g_{\varepsilon}}.\end{eqnarray*} Hence
\begin{equation}\int_{\widetilde{U}\bigcup A(r,2r)
}u^{2}dvol_{g_{\varepsilon}}\leq
\frac{\Lambda-R_{0}}{R_{1}-R_{0}}\int_{M
}u^{2}dvol_{g_{\varepsilon}}.
\end{equation} We have \begin{equation} \end{equation}
\begin{eqnarray*}\|v\|^{2}_{L^{2}(N)} &  = & \|\eta u\|^{2}_{L^{2}(N)} \geq  \| u\|^{2}_{L^{2}(N\backslash U(2r))}
\\  & \geq & (1-\frac{\Lambda-R_{0}}{R_{1}-R_{0}}) \| u\|^{2}_{L^{2}(M)}
\\  & \geq &
 \frac{R_{1}-\Lambda}{R_{1}-R_{0}} \|
 u\|^{2}_{L^{2}(M)}.\end{eqnarray*} For a $\rho\in [r,
 (2r)^{\frac{1}{11}}]$, set $\widehat{U}_{\rho}=\widetilde{U}\bigcup
 A(r, \rho)$.  Since $R_{g_{\varepsilon}}\geq R_{1}$ on
 $\widehat{U}_{\rho}$, we have \begin{eqnarray*}\lambda_{M}(g_{\varepsilon})\|
 u\|^{2}_{L^{2}(\widehat{U}_{\rho})}& = &
 4\int_{\widehat{U}_{\rho}}\langle \triangle u,
u\rangle dvol_{g_{\varepsilon}}+ \int_{\widehat{U}_{\rho}}
R_{g_{\varepsilon}}u^{2}dvol_{g_{\varepsilon}}\\ & =&
4\int_{\widehat{U}_{\rho}} |d
u|^{2}dvol_{g_{\varepsilon}}-4\int_{\partial
\widehat{U}_{\rho}}u\partial_{\nu}u d A +
\int_{\widehat{U}_{\rho}}
R_{g_{\varepsilon}}u^{2}dvol_{g_{\varepsilon}}\\ & \geq &
-4\int_{\partial \widehat{U}_{\rho}}u\partial_{\nu}u d A + R_{1}\|
 u\|^{2}_{L^{2}(\widehat{U}_{\rho})}.
\end{eqnarray*} Hence $$4\int_{\partial
\widehat{U}_{\rho}}u\partial_{\nu}u d A \geq
(R_{1}-\lambda_{M}(g_{\varepsilon}))\|
 u\|^{2}_{L^{2}(\widehat{U}_{\rho})}\geq
(R_{1}-\Lambda)\|
 u\|^{2}_{L^{2}(\widehat{U}_{\rho})}\geq 0.$$  By Lemma 3.7 in [BD2],
$$ \|u\|^{2}_{L^{2}(A(r, 2r))}\leq 10 r^{\frac{5}{2}}
\|u\|^{2}_{L^{2}(A(r, (2r)^{\frac{1}{11}}))}\leq 10
r^{\frac{5}{2}}\|u\|^{2}_{L^{2}(M)}.
$$ Thus \begin{equation}\frac{2}{r}\|u\|_{L^{2}(A(r, 2r))}\leq \zeta \|u\|_{L^{2}(M)}.
\end{equation} We have \begin{eqnarray*}\|d v\|^{2}_{L^{2}(N)}&=& \|d (\eta u)\|^{2}_{L^{2}(N)}\\ &\leq &
(\|\eta d u\|_{L^{2}(N)}+\|d \eta u\|_{L^{2}(N)})^{2}\\
&\leq & (\| d u\|_{L^{2}(M)}+\frac{2}{r}\| u\|_{L^{2}(A(r,
2r))})^{2}\\ &\leq & (\| d u\|_{L^{2}(M)}+\zeta\|
u\|_{L^{2}(M)})^{2}\\ &\leq & (1+\zeta)\| d
u\|^{2}_{L^{2}(M)}+\zeta(1+\zeta)\| u\|^{2}_{L^{2}(M)}
\end{eqnarray*} Since $R_{g_{\varepsilon}}>R_{1}>0$ on
$\widetilde{U}\bigcup U(2r)$, \begin{equation}
(R_{g}v,v)_{L^{2}(N)}\leq (R_{g_{\varepsilon}}\eta u,\eta
u)_{L^{2}(M)}\leq (R_{g_{\varepsilon}}u,u)_{L^{2}(M)}.
\end{equation} We obtain \begin{eqnarray*}& & 4\|d v\|^{2}_{L^{2}(N)}+(R_{g}v,v)_{L^{2}(N)}\\ &\leq & 4(1+\zeta)\| d
u\|^{2}_{L^{2}(M)}+\zeta(1+\zeta)\| u\|^{2}_{L^{2}(M)}+
(R_{g_{\varepsilon}}u,u)_{L^{2}(M)}\\ &\leq &
(1+\zeta)\lambda_{M}(g_{\varepsilon})\| u\|^{2}_{L^{2}(M)}-\zeta
R_{0}\| u\|^{2}_{L^{2}(M)} +\zeta(1+\zeta)\| u\|^{2}_{L^{2}(M)}\\
&\leq & [\lambda_{M}(g_{\varepsilon})+(\Lambda+1-R_{0})\zeta
+\zeta^{2}]\| u\|^{2}_{L^{2}(M)}.
\end{eqnarray*} Hence   we have \begin{eqnarray*}\lambda_{N}(g)& \leq
& \frac{\int_{N}(4|d
v|^{2}+R_{g}v^{2})dvol_{g}}{\int_{N}v^{2}dvol_{g}}\\ &\leq &
\frac{R_{1}-R_{0}}{R_{1}-\Lambda}[\lambda_{M}(g_{\varepsilon})+(\Lambda+1-R_{0})\zeta
+\zeta^{2}] \\ &=&
\lambda_{M}(g_{\varepsilon})+\frac{\Lambda-R_{0}}{R_{1}-\Lambda}
\lambda_{M}(g_{\varepsilon})+\frac{R_{1}-R_{0}}{R_{1}-\Lambda}[(\Lambda+1-R_{0})\zeta
+\zeta^{2}]\\ &\leq & \lambda_{M}(g_{\varepsilon}) +\varepsilon,
\end{eqnarray*} by (4.4), (4.5) and (4.14).  Thus
both Lemma 4.4 and Lemma 4.2 are proved. \end{proof}  \end{proof}
\vskip 3mm

\begin{lem} Let $N_{1}$ and $N_{2}$ be two compact oriented
4m-manifolds with $\overline{\lambda}_{N_{1}}\leq 0$,
$\overline{\lambda}_{N_{2}}\leq 0$. Let $M=N_{1}\sharp N_{2}$.
Assume that $\overline{\lambda}_{M}\leq 0$.  For any metrics
$g_{1}$ and $g_{2}$ on $N_{1}$ and $N_{2}$ respectively with
$\lambda_{N_{1}}(g_{1})=\lambda_{N_{2}}(g_{2})=-1$, and $0<
\varepsilon \ll 1$, there is a metric $g_{\varepsilon}$ on $M$
such that $$(1+\varepsilon)^{2m}
   (\overline{\lambda}_{N_{1}}(g_{1})^{2m}+
   \overline{\lambda}_{N_{2}}(g_{2})^{2m}+\varepsilon)\geq
   \overline{\lambda}_{M}(g_{\varepsilon})^{2m}.
$$ \end{lem}

\vskip 3mm

\begin{proof} Let $N=N_{1}\bigcup N_{2}$, and $p_{i}\in N_{i}$, $i=1, 2$.
Notations,  $r$, $\overline{r}$, $U(r)$,  $R_{0}$, $R_{1}$,
$\Lambda$, and $\zeta$, are  the same    as in the proof of Lemma
4.2. Here the only difference is that we use the set $\{p_{1},
p_{2}\}$ in stead of a point of $N$.  Denote  $U_{i}(r)=\{x\in
N_{i}| \text{dist}_{g_{i}}(x, p_{i})\leq r\}$. By Lemma 4.3,   for
any $0<\theta_{0}\ll 1$, for each $i$,  there is a metric
$\widetilde{g}_{i,\theta_{0}}$ on $A_{i}(r,
\frac{r}{2})=U_{i}(r)\backslash U_{i}(\frac{r}{2})$  satisfying
that $R_{\widetilde{g}_{i,\theta_{0}}}\geq R_{1}$,
 $\widetilde{g}_{i,\theta_{0}}$ agrees with $g_{i}$ near the boundary
$\partial U_{i}(r)$, and $\widetilde{g}_{i,\theta_{0}}$ agrees
with $dt^{2}+\delta(\theta_{0} )^{2}g_{0,
 1}$
 near the boundary  $\partial U_{i}(\frac{r}{2})\simeq S^{n-1}(1)$,
  where  $g_{0,
 1}$ is  the standard metric of sectional curvature 1
  on $S^{n-1}(1)$. From (4.6), we can choose  $\delta=\delta(\theta_{0}
)$ as   a function of $\theta_{0}$ in-dependent of $i$
 such that $\delta \ll \theta_{0}$.  Furthermore,  $$
|\text{Vol}_{\widetilde{g}_{i,\theta_{0}}}(A_{i}(r,
\frac{r}{2}))-\text{Vol}_{g_{i}}(U_{i}(r))|\longrightarrow 0,$$ if
$\theta_{0} \longrightarrow 0$. Note that $M$ is obtained by
gluing $N_{1}\backslash U_{1}(\frac{r}{2})$ and $N_{2}\backslash
U_{2}(\frac{r}{2})$ at $U_{1}(\frac{r}{2})\simeq
U_{2}(\frac{r}{2})$, i.e. $M=N_{1}\backslash U_{1}(\frac{r}{2})
\bigcup N_{2}\backslash U_{2}(\frac{r}{2})$. Define a metric
$g_{\theta_{0}}$ on $M$ by $g_{\theta_{0}}=g_{i}$ on
$N_{i}\backslash U_{i}(r)$, and
$g_{\theta_{0}}=\widetilde{g}_{i,\theta_{0}}$ on $A_{i}(r,
\frac{r}{2})$, which satisfies
$$|\text{Vol}_{g_{\theta_{0}}}(M)-\sum \text{Vol}_{g_{i}}(N_{i})|\leq
 \sum |\text{Vol}_{\widetilde{g}_{i,\theta_{0}}}(A_{i}(r,
\frac{r}{2}))-\text{Vol}_{g_{i}}(U_{i}(r))|\longrightarrow 0,$$
when $\theta_{0} \longrightarrow 0$. For any $\varepsilon >0$, by
letting $\theta_{0}\ll 1$ and $g_{\varepsilon}=g_{\theta_{0}}$, we
find  a metric $g_{\varepsilon}$  on $M$ with
$$|\text{Vol}_{g_{\varepsilon}}(M)-\sum \text{Vol}_{g_{i}}(N_{i})|\leq \varepsilon
.$$ By the same arguments as in the proof  of  Lemma 4.4, we have
$$-1-\varepsilon \leq
   \lambda_{M}(g_{\varepsilon}).$$ Thus $$(1+\varepsilon)^{2m}
   (\overline{\lambda}_{N_{1}}(g_{1})^{2m}+
   \overline{\lambda}_{N_{2}}(g_{2})^{2m}+\varepsilon)\geq
   \overline{\lambda}_{M}(g_{\varepsilon})^{2m}.
$$  \end{proof}

\vskip 3mm

\begin{proof}[Proof of Proposition  4.1]
  First, we assume that there is a
metric $h$ on $X$ with positive scalar curvature. By Lemma  4.2,
for any metric $g$ on $N$ and $0< \varepsilon \ll 1$, there exists
a metric $g_{\varepsilon}$ on $M$ such that
$$\lambda_{N}(g)-\varepsilon \leq \lambda_{M}(g_{\varepsilon}),  \ \
\ \ {\rm and} \ \ \ \
|\text{Vol}_{g}(N)-\text{Vol}_{g_{\varepsilon}}(M)|\leq
\varepsilon.$$ Thus we have
$$(\lambda_{N}(g)-\varepsilon)(\text{Vol}_{g}(N)+ \varepsilon)^{\frac{2}{n}} \leq
\lambda_{M}(g_{\varepsilon})\text{Vol}_{g_{\varepsilon}}(M)^{\frac{2}{n}}\leq
\overline{\lambda}_{M}, \ \ \ {\rm or}$$
$$(\lambda_{N}(g)-\varepsilon)(\text{Vol}_{g}(N)-
\varepsilon)^{\frac{2}{n}} \leq
\lambda_{M}(g_{\varepsilon})\text{Vol}_{g_{\varepsilon}}(M)^{\frac{2}{n}}\leq
\overline{\lambda}_{M}.$$ By letting $\varepsilon \longrightarrow
0$,
$$\lambda_{N}(g)\text{Vol}_{g}(N)^{\frac{2}{n}}\leq
\overline{\lambda}_{M}.$$ Thus $$ \overline{\lambda}_{N}=
  \sup\limits_{g\in \mathcal{M}} \overline{\lambda}_{N}(g)\leq
\overline{\lambda}_{M},
  $$ where $\mathcal{M}$ is the set of Riemannian
  metrics on $N$. Hence, we obtain (4.1).

  Now we assume that $\overline{\lambda}_{N}\leq 0$, $\overline{\lambda}_{X}\leq 0$,
   $\overline{\lambda}_{M}\leq 0$, and $n=4m$. We can  choose any two  metrics $g_{1}$ and
   $g_{2}$ on $N$ and $X$ respectively with  $\overline{\lambda}_{N}(g_{1})< 0$
    and $\overline{\lambda}_{X}(g_{2})< 0$. After re-scaling  them, we can assume
    $ \lambda_{N}(g_{1})=\lambda_{X}(g_{2})=-1$. By Lemma 4.5, for any $\varepsilon >0$, there
     exists a metric $g_{\varepsilon}$ on $M$ such that $$-(1+\varepsilon)
   (\overline{\lambda}_{N}(g_{1})^{2m}+
   \overline{\lambda}_{X}(g_{2})^{2m}+\varepsilon)^{\frac{1}{2m}}\leq
   \overline{\lambda}_{M}(g_{\varepsilon})\leq \overline{\lambda}_{M}.
$$ By letting $\varepsilon \longrightarrow 0$, we obtain $$-
   (\overline{\lambda}_{N}(g_{1})^{2m}+
   \overline{\lambda}_{X}(g_{2})^{2m})^{\frac{1}{2m}}\leq \overline{\lambda}_{M}.
$$  Thus $$-
   (\overline{\lambda}_{N}^{2m}+
   \overline{\lambda}_{X}^{2m})^{\frac{1}{2m}}=-
   ((\sup\limits_{g_{1}\in \mathcal{M}_{1}} \overline{\lambda}_{N}(g_{1}))^{2m}+
   (\sup\limits_{g_{2}\in \mathcal{M}_{2}}  \overline{\lambda}_{X}(g_{2}))^{2m})^{\frac{1}{2m}}\leq \overline{\lambda}_{M}.
$$  \end{proof}

\vskip 3mm

\begin{proof}[Proof of Theorem  1.4]
   By Theorem 1 in
   [Ta], the Spin$^{c}$ structure induced by a compatible almost complex
   structure on $(N, \omega)$ has Seiberg-Witten invariant equal to $\pm
   1$. By Lemma 1 in Section 3 of  [Le4], for any metric $g$ on $M$,  we can choose a Spin$^{c}$ structure on
   $M$ with non-vanishing Seiberg-Witten invariant, and
   $$[c_{1}^{+}]^{2}[M]\geq 2\chi(N)+3\tau(N)>0.$$ Thus, by Proposition
   3.2, we obtain $$\overline{\lambda}_{M}(g)\leq - \sqrt{
32\pi^{2}[c_{1}^{+}]^{2}[M]}\leq - \sqrt{
32\pi^{2}(2\chi(N)+3\tau(N))}.$$ Hence
$$\overline{\lambda}_{M}\leq - \sqrt{
32\pi^{2}(2\chi(N)+3\tau(N))}.$$ If $N$ admits a
K\"{a}hler-Einstein metric, then, by Proposition 4.1,
$$
-\sqrt{ 32\pi^{2}(2\chi(N)+3\tau(N))}=\overline{\lambda}_{N}\leq
\overline{\lambda}_{M}\leq -\sqrt{ 32\pi^{2}(2\chi(N)+3\tau(N))}.
 $$ Hence we obtain the conclusion.   \end{proof}

\vskip 3mm

\begin{lem} Let $M$ be a spin manifold with non-vanishing
$\widehat{A}$-genus, i.e. $\widehat{A}(M)\neq 0$. Then
$$\overline{\lambda}_{M}\leq 0.$$

\end{lem}

\vskip 3mm

\begin{proof} Since $\widehat{A}(M)\neq 0$,  for any metric $g'$ on $M$, there is a
 non-vanishing harmonic
 spinor $\phi \in \Gamma (S)$ with $\int_{M}|\phi|^{2}dvol_{g'}=1$,  where $S$ is the spin bundle.  The Bochner formula implies
 that $$0=\mathcal{D}^{2}\phi=\nabla^{*}\nabla \phi + \frac{R_{g'}}{4}\phi,$$
where $\mathcal{D}: \Gamma (S)\longrightarrow \Gamma (S)$ is the
Dirac operator.
 By (3.2),
\begin{eqnarray*}\lambda_{M}(g')& \leq &
\int_{M}(4|\nabla|\phi|_{\varepsilon}|^{2}+R_{g'}|\phi|_{\varepsilon}^{2})dvol_{g'}\\
& \leq & \int_{M}(4|\nabla \phi|^{2}+R_{g'}|\phi|^{2})dvol_{g'}
+\varepsilon^{2}\int_{M}R_{g'}dvol_{g'}=\varepsilon^{2}\int_{M}R_{g'}dvol_{g'}.\end{eqnarray*}
By letting $\varepsilon\longrightarrow 0$, we obtain that
$$\lambda_{M}(g')\leq 0, \ \ \ \ {\rm and} \ \ \ \
\overline{\lambda}_{M}\leq 0.$$
\end{proof}

\vskip 3mm

\begin{proof}[Proof of Theorem  1.5]
By Theorem 1 in
   [Ta], for any $i$,  the Spin$^{c}$ structure induced by a compatible almost complex
   structure on $(N_{i}, \omega_{i})$ has Seiberg-Witten invariant equal to $\pm
   1$. By Corollary 11 in [IL], for any metric $g$ on $M$,  there   is a monopole class $\alpha$ of $M$ satisfying that
   $$[\alpha^{+}]^{2}[M]\geq \sum_{i=1}^{\ell}c_{1}^{2}[N_{i}].$$ Thus, by Proposition
3.2, we obtain $$\overline{\lambda}_{M}(g)\leq - \sqrt{
32\pi^{2}[\alpha^{+}]^{2}[M]}\leq - \sqrt{
32\pi^{2}\sum_{i=1}^{\ell}c_{1}^{2}[N_{i}]}.$$ Hence
$$\overline{\lambda}_{M}\leq - \sqrt{
32\pi^{2}\sum_{i=1}^{\ell}c_{1}^{2}[N_{i}]}.$$

Now, we assume that, for each $i$,  $N_{i}$ admits a
K\"{a}hler-Einstein metric $g_{i}$. If the scalar curvature of
$g_{i}$ is negative, we have already known that
$\overline{\lambda}_{N_{i}}=\overline{\lambda}_{N_{i}}(g_{i})=-
\sqrt{ 32\pi^{2}c_{1}^{2}(N_{i})}$ by Theorem 1.1. If the scalar
curvature of $g_{i}$ is zero, then $N_{i}$ is a $K3$-surface from
the hypothesis (cf. [BHPV]). By Lemma 4.6, we obtain
$0=\overline{\lambda}_{N_{i}}(g_{i})\leq
\overline{\lambda}_{N_{i}}\leq 0$.  Thus
$\overline{\lambda}_{N_{i}}=0=- \sqrt{
32\pi^{2}c_{1}^{2}(N_{i})}$. By Proposition 4.1, $$- \sqrt{
32\pi^{2}\sum_{i=1}^{\ell}c_{1}^{2}[N_{i}]}=-\sqrt{\sum_{i=1}^{l}\overline{\lambda}_{N_{i}}^{2}}\leq
\overline{\lambda}_{M}\leq - \sqrt{
32\pi^{2}\sum_{i=1}^{l}c_{1}^{2}[N_{i}]}.
$$ This proves the desired result.    \end{proof}

\vskip 3mm

\begin{proof}[Proof of Corollary   1.6]
  Note that $N_{i}$  are spin manifolds (See [J]), and
 thus is
 $M$. By Lemma 4.6, $\overline{\lambda}_{M}\leq 0$ as $\widehat{A}(M)\neq 0$.
    Since  $X_{1} \cdots X_{l}$ are
simply connected  compact oriented spin n-manifolds with
$\widehat{A}(X_{j})=0$, $n\geq 8$,
  and $n=0{\rm mod}4$, for any $X_{j}$, there is a metric $h_{j}$ on
   $X_{j}$ with positive scalar curvature from Theorem A in [St].
  Note that $(N_{i}, g_{i})$ are
  Ricci-flat Einstein manifolds with $\widehat{A}(N_{i})\neq 0$ (See
  [J]). By Lemma 4.6, $0=\overline{\lambda}_{N_{i}}(g_{i})\leq
  \overline{\lambda}_{N_{i}}\leq 0$.  By Proposition
4.1,
$$0 \leq  \overline{\lambda}_{\sharp_{i=1}^{l_{1}}N_{i}}\leq\overline{\lambda}_{M}\leq 0.$$
  We obtain the conclusion. \end{proof}

\vskip 3mm

\section{Proof of Proposition  1.7}

\vskip 3mm

\begin{proof} [Proof of Proposition  1.7]  If it is  not true, there
exists a sequence of metrics $\{g_{k}\}\subset
\mathcal{M}_{(\Lambda, D)}$ such that
$$- \sqrt{32\pi^{2}
(2\chi(M)+3\tau(M))}-\frac{1}{k}
\leq\overline{\lambda}_{M}(g_{k})$$ but $g_{k}$ can  never  be
deformed to a complex hyperbolic metric through the Ricci flow for
every $k$.

Since $\chi(M)>0$, there is a positive constant $v$ independent of
$k$ such that $\text{Vol}_{g_{k}}(M)\geq v$ by the
Gauss-Bonnett-Chern theorem. By the Cheeger-Gromov theorem (cf
[A]),
 $\{g_{k}\}$ has a  $C^{1,\alpha}$-convergence  subsequence, denoted by $\{g_{k}\}$ also.
Therefore, there are diffeomorphisms $F_{k}$ of $M$ such that
 a subsequence of
 $\{F_{k}^{*}g_{k}\}$ converges, in the $C^{1,\alpha}$-topology on
 $M$, to a $C^{1,\alpha}$-metric $g_{\infty}$.  In fact,
 $\{F_{k}^{*}g_{k}\}$ converges in the $L^{2,p}$-topology, for any $p\geq 1$, and
 $g_{\infty}$ is a $L^{2,p}$-metric (See [A] for details). Thus
 $\overline{\lambda}_{M}(g_{\infty})$ is well defined satisfying
 that
 $$
-\sqrt{32\pi^{2} (2\chi(M)+3\tau(M))}
\leq\overline{\lambda}_{M}(g_{\infty}).$$ This together with
Theorem 1.1 implies that
$$\overline{\lambda}_{M}(g_{\infty})=-
\sqrt{32\pi^{2} (2\chi(M)+3\tau(M))},$$ and $g_{\infty}$ is a
K\"{a}hler-Einstein metric with negative scalar curvature. By
Theorem 5 in [Le1] $\chi(M)\ge 3\tau (M)$. This together with the
assumption $\chi (M) \in [\frac 32 \tau (M), 3\tau (M)]$ implies
that $\chi(M)=3\tau (M)$. By Theorem 5 in [Le1] once again we know
that $ g_{\infty}$ is a complex hyperbolic metric.

To prove the metric can be deformed to a complex hyperbolic metric
through the Ricci flow, we need to smooth the
$C^{1,\alpha}$-convergence to a $C^2$-convergence by Ricci flow.
By the main theorem in [BOR] (See Theorem 5.1 in [Fu] for this
version), given any $1\gg \epsilon
>0$ and $j\in \mathbb{N}$, there exists a constant $C(j, \epsilon)$
and a smoothing operator $S_{\epsilon}: \mathcal{M}_{(\Lambda,
D)}\longrightarrow \mathcal{M}_{(2\Lambda, 2D)}$ such that
\begin{enumerate}
 \item  $ \|S_{\epsilon}(g)-g\|_{C^{0}}<\epsilon$,
 \item   $ \|\nabla^{S_{\epsilon}(g)}-\nabla^{g}\|_{C^{0}}<\epsilon$,
\item  $\|\nabla^{j}{\rm Rm}(S_{\epsilon}(g))\|_{C^{0}}<C(j,
\epsilon)\|{\rm Rm}(g)\|_{C^{0}}$,
 \end{enumerate} where ${\rm Rm}(g)$ is the curvature operator of $g$.
The proof of this result is by  considering the Ricci-flow
evolution equation
 with initial metric $g\in  \mathcal{M}_{(\Lambda, D
 )}$  $$\begin{array}{ccc}\frac{\partial}{\partial t}g(t)=-2{\rm Ric}(g(t)) \\
   g(0)=g,
\end{array}  $$ and letting  $S_{\epsilon}(g)=g(\epsilon)$. By
using the operator $S_{\epsilon}$ to metrics $g_{k}$, we obtain a
sequence of metrics $\{S_{\epsilon}(g_{k})\}\subset
\mathcal{M}_{(2\Lambda, 2D)}$.  Let
$\widetilde{g}_{k}=S_{\epsilon}(g_{k})$.  By the claim in $\S$2.3
of [Pe1],  $\overline{\lambda}_{M}(g)$ is non-decreasing along the
Ricci flow if $\overline{\lambda}_{M}(g)\leq 0$. Thus $$-
\sqrt{32\pi^{2} (2\chi(M)+3\tau(M))}-\frac{1}{k}
\leq\overline{\lambda}_{M}(g_{k})\leq
\overline{\lambda}_{M}(\widetilde{g}_{k}).$$

By the Cheeger-Gromov Theorem again, there are diffeomorphisms
$\widetilde{F}_{k}$ of $M$ such that
 a subsequence of $\{\widetilde{F}_{k}^{*}\widetilde{g}_{k}\}$, saying
$\{\widetilde{F}_{k}^{*}\widetilde{g}_{k}\}$ again, which
converges in the $C^{1,\alpha}$-topology in $M$ to a
$C^{1,\alpha}$-metric $\widetilde{g}_{\infty}$, a
K\"{a}hler-Einstein metric with negative scalar curvature by
Theorem 1.1. Since $\|\nabla {\rm
Rm}(\widetilde{g}_{k})\|_{C^{0}}<C(1, \epsilon)\Lambda$, by the
Arzela-Ascoli Theorem we get a sub-sequence of $\{{\rm
Rm}(\widetilde{F}_{k}^{*}\widetilde{g}_{k})\}$ which
$C^{0}$-converges to ${\rm Rm}(\widetilde{g}_{\infty})$.
Therefore, $\{\widetilde{F}_{k}^{*}\widetilde{g}_{k}\}$
$C^{2}$-converges to $\widetilde{g}_{\infty}$.  As above by [Le1]
$ \widetilde{g}_{\infty}$ is a complex hyperbolic metric. Note
that the sectional curvature $K(\widetilde{g}_{\infty})$ of a
complex hyperbolic metric is negative, i.e. there are constants
$\mu_{1}$ $\mu_{2}$ such that $-\mu_{1}^{2}\leq
K(\widetilde{g}_{\infty})\leq -\mu_{2}^{2}$. Thus, for $k\gg 1$,
we have $-2\mu_{1}^{2}\leq K(\widetilde{g}_{k})\leq
-\frac{1}{2}\mu_{2}^{2}$. Moreover, the Einstein tensors satisfy
that
$$T_{\widetilde{g}_{k}}={\rm Ric}(\widetilde{g}_{k})-\frac{R_{\widetilde{g}_{k}}}{4}\widetilde{g}_{k}\longrightarrow
0$$ in the $C^{0}$-sense  when $k\longrightarrow \infty$. By the
corollary of Theorem 1.1 in [Ye], for a $k\gg 1$,
$\widetilde{g}_{k}$ can be deformed to an Einstein metric, which
is complex hyperbolic metric by [Le1] again.

Note that we first deform $g_{k}$ to $\widetilde{g}_{k}$ through
the Ricci flow, then deform  $\widetilde{g}_{k}$ to a complex
hyperbolic metric  through the Ricci flow again. A contradiction.
The desired result follows.  \end{proof}

\vskip 14mm

\vskip 20mm

%%%%%%%%%%%%%%%%%%%%%%%%%%%

\end{document}